\documentclass[12pt,leqno,draft]{article}

\newtheorem{theorem}{Theorem}
\newtheorem{lemma}[theorem]{Lemma}
\newtheorem{proposition}[theorem]{Proposition}
\newtheorem{sublemma}[theorem]{Sublemma}
\newtheorem{definition}[theorem]{Definition}
\newtheorem{corollary}[theorem]{Corollary}
\newtheorem{problem}[theorem]{Problem}
\newtheorem{remark}[theorem]{Remark}
\newtheorem{claim}[theorem]{Claim}
\newtheorem{assumptions}[theorem]{Assumptions}
\newtheorem{examples}[theorem]{Examples}
\newtheorem{question}[theorem]{Question}
\newtheorem{sassumptions}[theorem]{Standing Assumptions}
\newtheorem{sassumption}[theorem]{Standing Assumption}
\newtheorem{conjecture}[theorem]{Conjecture}

\newcommand{\begintheorem}{\addtocounter{equation}{1}\begin{theorem}}
\newcommand{\beginlemma}{\addtocounter{equation}{1}\begin{lemma}}
\newcommand{\beginproposition}{\addtocounter{equation}{1}\begin{proposition}}
\newcommand{\beginsublemma}{\addtocounter{equation}{1}\begin{sublemma}}
\newcommand{\begindefinition}{\addtocounter{equation}{1}\begin{definition}}
\newcommand{\begincorollary}{\addtocounter{equation}{1}\begin{corollary}}
\newcommand{\beginproblem}{\addtocounter{equation}{1}\begin{problem}}
\newcommand{\beginremark}{\addtocounter{equation}{1}\begin{remark}}
\newcommand{\beginclaim}{\addtocounter{equation}{1}\begin{claim}}
\newcommand{\beginassumptions}{\addtocounter{equation}{1}\begin{assumptions}}
\newcommand{\beginexamples}{\addtocounter{equation}{1}\begin{examples}}
\newcommand{\beginquestion}{\addtocounter{equation}{1}\begin{question}}
\newcommand{\beginsassumptions}{\addtocounter{equation}{1}\begin{sassumptions}}
\newcommand{\beginsassumption}{\addtocounter{equation}{1}\begin{sassumption}}
\newcommand{\beginconjecture}{\addtocounter{equation}{1}\begin{conjecture}}

\begin{document}

\title{Some remarks about homeomorphisms, ``energy'', and so on}

\author{S. Semmes}

\date{}

\maketitle

	John Ball keeps asking questions of the following sort.
Suppose that one has a homeomorphism from a domain in a Euclidean
space onto its image in the same Euclidean space.  Assume also that
the homeomorphism has ``finite energy'' with respect to some
reasonable functional, which would normally entail something like
first distributional derivatives in $L^p$ for some $p$, $1 \le p <
\infty$, and the inverses of the associated differentials lying in some
$L^q$, $1 \le q < \infty$.  The latter might be controlled in terms
of integrability properties of the reciprocal of the Jacobian of
the mapping (the determinant of the differential).  Under conditions
such as these, what kinds of approximations of the homeomorphism
can one make by more regular homeomorphisms, approximations which
respect similar integrability conditions for the differentials and
their inverses?

	Let us restrict ourselves to dimensions $3$ and lower, since
all sorts of strange things happen in dimensions $4$ and larger, and
since dimensions less than or equal to $3$ are physically relevant
(elasticity theory, etc., as Dr.\ Ball well knows).  This type of
issue in dimension $1$ can be treated in a direct manner, by writing
the mapping as the integral of its derivative, and so we focus on
dimensions $2$ and $3$.

	There are very famous results about approximating
homeomorphisms by piecewise-linear homeomorphisms in dimensions $2$
and $3$.  See \cite{Bing, Moise}.  More precisely, these are
approximations in $C^0$ senses, which are already quite nontrivial and
useful in the study of topology.  The results include \emph{relative}
versions, in which a homeomorphism is regularized in some parts while
not changing it on other parts where it is already regular.

	What about approximations which also respect the differential
in some manner?

	A basic strategy in making approximations of a function in a
Sobolev space is to first choose a set on which the function behaves
nicely, and whose complement has small measure.  The nice behavior
might involve boundedness of the first derivatives, continuity of the
first derivatives, and so forth.  Although the complement has small
measure, it typically need not have any simple structure, but could be
quite scatterred and irregular.  The second step would be to modify
the function on the small set.  If there are no additional constraints
on the function, then there are relatively simple techniques of
extension and regularization.  However, if one wants the result to be
a homeomorphism, then this second step becomes much more complicated.

	For more on approximations of homeomorphisms, if not exactly
of this form, see \cite{D-S, Luu3, Sullivan, TV2, V1}.

	This type of conundrum seems a bit odd to me, in that for a
number of basic situations that would arise in a simple way, there
should not be as much trouble.  One way to look at this is that often
there is something like a one-parameter family of mappings.  A related
issue is that in topology a basic point is often to construct
\emph{isotopies} between homeomorphisms, i.e., a continuous family of
homeomorphisms.  This is quite different from \emph{homotopies}, which
are continuous families of mappings which are not required to be
injective, even if the mappings that the homotopy goes between are
injective.

	In particular, there are a number of results to the effect
that homeomorphisms which are close in a $C^0$ sense can be connected
by an isotopy, as in \cite{Cer, E-K, Fisher, Hamstrom, Kister1,
Kister2, Luu1, Sullivan}, and where the isotopy stays close to the
original homeomorphisms.  This is part of the reason that a $C^0$
approximation can be useful, since otherwise it seems rather weak.

	In another direction, let us recall a famous result of Hatcher
\cite{Hatcher1, Hatcher2, Lau} concerning embedded two-dimensional
spheres in ${\bf R}^3$.  More precisely, one considers smoothly embedded
two-dimensional spheres in ${\bf R}^3$, and the space of these can be
locally identified with the space of smooth real-valued functions on the
$2$-sphere, because of the tubular neighborhood theorem.  Hatcher's
result says that this space is \emph{contractable}.

	For the analogous question in the plane, the Riemann mapping
theorem can be used.  Paul Schweitzer keeps asking about possible 
analytic proofs for two-dimensional spheres in ${\bf R}^3$.

	Note that embedded $2$-spheres in ${\bf R}^3$ can be viewed as
the boundaries of solid $3$-dimensional balls, just as Jordan curves
in the plane can be viewed as boundaries of $2$-dimensional disks.

	For another very interesting direction along the lines of
analysis and complexity of shapes, see \cite{C-K-S, Ku-S1, Ku-S2}.

\end{document}